\documentclass[12pt]{article}

\usepackage{amsmath}
\usepackage{theorem}
\usepackage{amssymb}
\usepackage{latexsym}

\usepackage{a4wide}

\newtheorem{thm}{Theorem}
\newtheorem{lem}[thm]{Lemma}

{\theorembodyfont{\rmfamily} }
\newenvironment{rem}{\noindent{\bf Remark.}}{\newline}
\newenvironment{pf}{\noindent{\bf Proof.}}{\hbox{}\hfill $\Box$}

\newcommand{\Q}{\mathbb{Q}}

\newcommand{\F}{\mathbb{F}}

\newcommand{\ad}{\mathrm{\mathop{ad}}}
\newcommand{\Tr}{\mathrm{\mathop{Tr}}}
\newcommand{\Aut}{\mathrm{\mathop{Aut}}}
\newcommand{\GL}{\mathrm{\mathop{GL}}}
\newcommand{\Der}{\mathrm{\mathop{Der}}}
\newcommand{\Inn}{\mathrm{\mathop{Inn}}}

\begin{document}

\title{Classification of solvable Lie algebras}
\author{W. A. de Graaf, \\
School of Mathematics and Statistics,\\
University of Sydney, \\
Australia\\
email: {\tt wdg@maths.usyd.edu.au}}

\date{}
\maketitle

\begin{abstract}
In this paper we illustrate some simple ideas that can be used for 
obtaining a classification of small-dimensional solvable Lie algebras.
Using these we obtain the classification of $3$ and $4$ dimensional
solvable Lie algebras (over fields of any characteristic). Precise
conditions for isomorphism are given. 
\end{abstract}

Solvable Lie algebras have been classified by G. M. Mubarakzjanov (upto dimension
$6$ over a real field, \cite{mub1}, \cite{mub2}, see also \cite{pwz}), and
by J. Patera and H. Zassenhaus (upto dimension $4$ over any perfect field, 
\cite{patzas2}). In this paper we explore the possibility of using the computer
to obtain a classification of solvable Lie algebras. The possible advantages of 
this are clear. The problem of classifying Lie algebras needs a systematic 
approach, and the more the computer is invloved, the more systematic the
methods have to be. However, the drawback is that the computer can only handle
finite data. For example, below we will consider orbits of the action of 
the automorphism group of a Lie algebra on the algebra of its derivations.
Now, if the ground field is infinite, then I know of no algorithm for
obtaining these orbits. In our approach we use the computer (specifically the
technique of Gr\"{o}bner bases) to decide isomorphism of Lie algebras, and to
obtain explicit isomorphisms if they exist. \par
The procedure that we use to classify solvable Lie algebras is based on some
simple ideas, which are described in Section \ref{sec1} (and for which we do
not claim any originality). Then in Section \ref{sec2} we describe the use of
Gr\"{o}bner bases. In Section \ref{sec3} solvable Lie algebras of dimension
$3$ are classified. In Section \ref{sec4} the same is done for dimension $4$.
We show that our classification in dimension $4$ differs slightly from the
one found in \cite{patzas2} (i.e., we find a few more Lie algebras). \par
For doing the explicit calculations reported here we have used the computer 
algebra system {\sf Magma} (\cite{magma}).\par

\section{General lemmas}\label{sec1}

In the sequel we denote the field we work with by $F$. 

\begin{lem}
Let $L$ be a solvable Lie algebra. Then there is a subalgebra $K\subset L$
of codimension $1$, and a derivation $d$ of $K$ such that 
$L= Fx_d\oplus K$. Here $[x_d,y]=d(y)$ for all $y\in K$. Moreover, if 
$L$ is not Abelian, then $d$ and $K$ 
can be chosen such that $d$ is an outer derivation of $K$.
\end{lem}

\begin{pf}
Let $K$ be any subspace of codimension $1$ containing $[L,L]$, and 
let $x\in L$ span a complementary subspace. Then $K$ is an ideal of $L$
and we get the result with $d=\ad x$ and $x_d=x$.\par
The proof of the second statement is by induction on $\dim L$. If $\dim L=2$
then the statement is clear. Now suppose $\dim L = n>2$ and write $L = 
Fy\oplus K$. Suppose that $\ad y$ is an inner derivation of $K$, i.e., 
$\ad y =\ad u$ for some $u\in K$. Set $z=y-u$; then $L= K\oplus Fz$ and 
$[z,K]=0$. So $K$ is non-Abelian and by induction we have $K= Fx\oplus K_1$, 
where $\ad x|_{K_1}$ is an outer derivation. Set $K_2 = K_1\oplus Fz$, then
also $\ad x|_{K_2}$ is an outer derivation, and $L = Fx\oplus K_2$.
\end{pf}

\begin{lem}\label{lem2}
Let $K$ be a solvable Lie algebra and $d_1,d_2$ derivations of $K$.
Set $L_i=Fx_{d_i}\oplus K$, $i=1,2$. Suppose that there is an automorphism 
$\sigma$ of $K$ such that $\sigma d_1\sigma^{-1} = \lambda d_2$, for some
scalar $\lambda\neq 0$. Then $L_1$ and $L_2$ are isomorphic.
\end{lem}

\begin{pf}
Define a linear map $\tilde{\sigma} : L_1\to L_2$ by $\tilde{\sigma}(y)
=\sigma(y)$ for $y\in K$ and $\tilde{\sigma}(x_{d_1}) = \lambda x_{d_2}$. Then 
$\tilde{\sigma}$ is a bijective linear map. The fact that it is an isomorphism
can be established by direct verification.
\end{pf}

The procedure based on these lemmas is as follows. Let $K$ be a 
solvable Lie algebra of dimension $n$. We compute the automorphism group 
$\Aut(K)$ of $K$ and the derivation algebra $\Der(K)$ of $K$. We denote
the subalgebra of inner derivations by $\Inn(K)$. It is straightforward 
to see that Lie algebras defined by derivations in the same coset of 
$\Inn(K)$ in $\Der(K)$ are isomorphic. Now the group $G(K)=F^*\times
\Aut(K)$ acts on
the cosets $d+\Inn(K)$ for $d\in \Der(K)$ by $(\lambda,\sigma) \cdot d 
+\Inn(K)= \lambda \sigma d \sigma^{-1}+\Inn(K)$. We compute orbit 
representatives of the action of $G(K)$ on $\Der(K)/\Inn(K)$. For every 
such representative we get a solvable Lie algebra of dimension $n+1$.
Subsequently we weed out the isomorphic ones. \par
When doing this we often deal with Lie algebras given by a multiplication
table containing parameters. An easy trick that often works to reduce
the number of parameters is to consider a diagonal base change.
Let $\{x_1,\ldots,x_n\}$ be a basis of $L$, and set $y_i =\alpha_i x_i$.
Then write down the multiplication table of $L$ with respect to the $y_i$.
Often it is possible to choose the $\alpha_i$ in such a way that we can get
rid of one or more parameters.\par
When $K$ is Abelian of dimension $n$ we have that $\Der(K) = M_n(F)$
and $\Aut(K) = \GL(n,F)$. In this case representatives of the orbits
of $\Aut(K)$ are known, by the following well-known theorem (for a proof
we refer to \cite{harhaw}). 
\begin{thm}\label{canformthm}
Let $A$ be an $n\times n$-matrix over a field $F$. Then $A$ is similar
over $F$ to a unique block-diagonal matrix, containing the blocks $C(f_1),
\ldots, C(f_s)$ where $C(f_k)$ is the companion matrix of the
non-constant monic polynomial $f_k$, and $f_k | f_{k+1}$ for $1\leq k\leq s-1$.
\end{thm}
The unique block-diagonal matrix is called the {\em rational canonical 
form} of $A$. \par
We end this section with some remarks on notational conventions that we use.
In this paper we usually describe an $n$-dimensional Lie algebra 
by giving its multiplication table with respect to a basis, which on most
occasions is denoted $x_1,\ldots, x_n$. In these multiplication tables we 
use the convention that products which are not listed are zero. Also
when representing a linear map by a matrix we always use the column 
convention.

\section{Constructing isomorphisms}\label{sec2}

One of the main problems when classifying Lie algebras is to decide
whether two of them are isomorphic. A very convenient tool for doing that
is Gr\"{o}bner bases (cf. \cite{gela}). (For an introduction into Gr\"{o}bner
bases we refer to \cite{clo}.) By way of example we describe how this works. \par
Consider the $3$-dimensional Lie algebra $L_1$ with basis $x_1,x_2,x_3$ and
multiplication table
$$[x_1,x_2]=x_2,~ [x_1,x_3]=ax_3,$$
and the $3$-dimensional Lie algebra $L_2$ with basis $y_1,y_2,y_3$ and
multiplication table
$$[y_3,y_1]=y_2,~ [y_3,y_2] = by_1+y_2.$$
The question is whether $L_1$ and $L_2$ are isomorphic, and if 
so for which values of $a,b$. In that case we would also like to have an 
explicit isomorphism. An isomorphism will map the nilradical of 
$L_1$ onto the nilradical of $L_2$. So an isomorphism $\phi : L_1\to
L_2$ has the form $\phi(x_1)=a_1y_1+a_2y_2+a_3y_3$, 
$\phi(x_2)=b_1y_1+b_2y_2$, $\phi(x_3)=c_1y_1+c_2y_2$. Now this is an
isomorphism if and only if the following polynomial equations are satisfied
$$ba_3b_2-b_1=0, ~ a_3b_1+a_3b_2-b_2=0,~ ba_3c_2-ac_1=0, ~ 
a_3c_1+a_3c_2-ac_2=0,$$
and
$$D_1a_3-1=0,~ D_2(b_1c_2-b_2c_1)-1 =0.$$
The last two equations are added to ensure that the determinant is nonzero.
Now in {\sf Magma} we compute a Gr\"{o}bner basis of the ideal of
$\Q[D_1,D_2,a_1,a_2,a_3,b_1,b_2,c_1,c_2,a,b]$ generated
by the left hand sides of these equations. We use the lexicographical 
ordering. This leads to a Gr\"{o}bner basis with a triangular structure,
which on many occasions makes it possible to find an explicit solution.
Also, we let $a,b$ be the smallest variables in the ordering; this 
makes it likely that the Gr\"{o}bner basis contains polynomials in only
$a$ and $b$ (cf. \cite{clo}, Chapter 3, Theorem 2).
From these we can derive necessary conditions for isomorphism. 
So in {\sf Magma} we do the following
\begin{verbatim}
> P<D1,D2,a1,a2,a3,b1,b2,c1,c2,a,b>:= PolynomialRing( Rationals(), 11 ); 
> r:= [ b*a3*b2-b1, a3*b1+a3*b2-b2, b*a3*c2-a*c1, a3*c1+a3*c2-a*c2, 
> D1*a3-1, D2*(b1*c2-b2*c1)-1 ];  
> I:= Ideal( r );
> g:= GroebnerBasis( I );
> g;
[
    D1 - a*b - b - 1,
    D2*b2*c2*a - D2*b2*c2 + a + 1,
    D2*b2*c2*b + 1/4*D2*b2*c2 - 1/2*a*b - 1/2*b - 1/4,
    a3 - a - 1,
    b1 - b2*a*b - b2*b,
    c1 + c2*a*b + c2*b + c2,
    a^2*b + 2*a*b + a + b
]
\end{verbatim}
From the last equation we get
\begin{equation}\label{eqn1}
b=-\frac{a}{(a+1)^2}.
\end{equation}
From this we also see that the algebras are not isomorphic if $a=-1$.
Solutions to the fourth, fifth and sixth equations are easily found, e.g., 
$a_1=a_2=0$, $a_3=a+1$, $b_1=ab+b$, $b_2=1$, $c_1=ab+b+1$, $c_2=-1$. 
By direct verification we get that this indeed defines an isomorphism,
if (\ref{eqn1}) holds, and $a\neq 1$ (otherwise the determinant is zero). 
If $a=1$ then by a separate calculation we get that the Gr\"{o}bner basis
is $\{1\}$. So in that case $L_1$ and $L_2$ are not isomorphic. 
The conclusion is that $L_1\cong L_2$ if and only if (\ref{eqn1}) and
$a\neq 1$. Moreover, in that case we also have an explicit isomorphism.\par
In the above discussion we have taken the ground field to be $\Q$. However,
the conclusion holds over any field of characteristic $0$, since over
any such field the Gr\"{o}bner basis will be the same. We can also easily
reach the same conclusion for any field of characteristic $p>0$. For that
we note that the input polynomials $g_i$ are defined over any field. Now
by using the {\sf Magma} function {\sf Coordinates} we can find polynomials 
$p_i$ such that $\sum p_ig_i=a^2b + 2ab + a + b$. The coefficients
of the $p_i$ are rational numbers. So from these coefficients we find a 
finite set of characteristics over which they are not defined. We then
have to do the computation separately over fields of those characteristics.
Over all other fields we have that the ideal generated
by the input polynomials contains $a^2b + 2ab + a + b$. Hence for all
those fields we reach the same conclusion (since the explicit isomorphism 
is defined over any field). In our example we have
\begin{verbatim}
> Coordinates( I, g[7] );
[
    -D2*a3*c1*a^2*b - D2*a3*c1*a*b - D2*a3*c1*a - D2*a3*c2*a^2*b^2 
    - D2*a3*c2*a^2*b - 2*D2*a3*c2*a*b^2 - 2*D2*a3*c2*a*b - D2*a3*c2*a - 
    D2*a3*c2*b^2 - D2*a3*c2*b + D2*c1*a^3*b + 2*D2*c1*a^2*b + D2*c1*a^2 + 
    D2*c1*a*b + D2*c2*a^3*b + D2*c2*a^2*b + D2*c2*a^2 - D2*c2*a*b - 
    D2*c2*a - D2*c2*b,
    -D2*a3*c1*a*b - D2*a3*c1*b - D2*a3*c1 - D2*a3*c2*a*b^2 - D2*a3*c2*a*b 
    - D2*a3*c2*b^2 - 2*D2*a3*c2*b - D2*a3*c2 + D2*c1*a^2*b + 2*D2*c1*a*b + 
    D2*c1*a + D2*c1*b + D2*c2*a^2*b + 2*D2*c2*a*b + D2*c2*a,
    D2*a3*b1*a*b + D2*a3*b1*b + D2*a3*b1 + D2*a3*b2*a^2*b^2 + 
    2*D2*a3*b2*a*b^2 + 2*D2*a3*b2*a*b + D2*a3*b2*b^2 + 2*D2*a3*b2*b + 
    D2*a3*b2 - D2*b1*a^2*b - 2*D2*b1*a*b - D2*b1*a - D2*b1*b,
    D2*a3*b1*a*b + D2*a3*b1*b + D2*a3*b1 + D2*a3*b2*a^2*b^2 + 
    D2*a3*b2*a*b^2 + 2*D2*a3*b2*a*b + D2*a3*b2*b + D2*a3*b2 - 
    D2*b1*a^2*b - 2*D2*b1*a*b - D2*b1*a - D2*b1*b - D2*b2*a*b - 
    2*D2*b2*b - D2*b2,
    0,
    -a^2*b - 2*a*b - a - b
]
\end{verbatim}

The coefficients of these polynomials are all integers. Hence 
the above conclusion holds for all fields.

\section{The $3$-dimensional case}\label{sec3}

There are only two (isomorphism classes of) Lie algebras of dimension $2$.\par
First we consider the Lie algebra $K$ spanned by $x_1,x_2$ with $[x_1,x_2]=0$.
Then $\Aut(K) = \GL(2,F)$, and $\Der(K) = M_2(F)$ (i.e., the space of all
$2\times 2$-matrices). In this case the rational canonical form of an
element in $\Der(K)$ is either
$$\lambda \begin{pmatrix} 1&0\\0&1\end{pmatrix},$$
or 
$$\begin{pmatrix} 0&a\\1&b\end{pmatrix}.$$
Fist of all, $\lambda=0$ 
gives a Lie algebra that is the three dimensional Abelian Lie algebra.
If $\lambda\neq 0$ then by
Lemma \ref{lem2} we may divide by $\lambda$ and get the Lie algebra $L^2$
spanned by $x_1,x_2,x_3$ and nontrivial brackets $[x_3,x_1]=x_1$, 
$[x_3,x_2]=x_2$. In the third case we get the Lie algebras 
$L_{a,b}$ spanned by $x_1,x_2,x_3$ and multiplication table 
$$[x_3,x_1] = x_2, ~[x_3,x_2]= ax_1+bx_2.$$ 
From a Gr\"{o}bner basis computation we get that $L_{a,b}\cong L_{c,d}$ 
implies $ad^2-b^2c=0$ and $c_3^2c-a=0$, for some nonzero $c_3\in F$. 
Furthermore, this holds over any field since the coordinates of these 
polynomials (with respect to the input basis) all have integer 
coefficients.\par
Also by applying a diagonal base change we see   
that $L_{a,b}\cong L_{\alpha^2a,\alpha b}$ for $\alpha\neq 0$
(the base change is $y_1=x_1$, $y_2=\alpha^{-1}x_2$
$y_3=\alpha^{-1} x_3$). Now we distinguish two cases:
\begin{enumerate}
\item $b\neq 0$. Then $L_{a,b}\cong L_{a,1}$. So we get the
class of Lie algebras $L_a^3 = L_{a,1}$. The above discussion implies that
$L_a^3\cong L_b^3$ if and only if $a=b$.
\item $b=0$. We get the class of Lie algebras $L_a^4= L_{a,0}$. In this case
$L_a\cong L_b$ if and only if $a=\alpha^2 b$ for some nonzero $\alpha\in F$.
\end{enumerate}
The Lie algebras $L^3_{a}$ and $L^4_b$ are never isomorphic. This can be 
established by a Gr\"{o}bner basis computation. It can also be shown in 
the following way. Suppose that $L^3_a\cong L^4_c$. Then $ad^2-b^2c=0$
amounts to $c=0$. However, $L^4_0$ is nilpotent, and $L^3_a$ is not.\par
The other $2$-dimensional Lie algebra (with basis $x_1,x_2$
and $[x_1,x_2]=x_2$) does not have to be considered, as it has no
outer derivations. \par
Summarising we get the following solvable Lie algebras of dimension $3$:
\begin{itemize}
\item[$L^1$] $[x_3,x_1]=0,~ [x_3,x_2]=0,~[x_1,x_2]=0.$
\item[$L^2$] $[x_3,x_1]=x_1,~ [x_3,x_2]=x_2,~[x_1,x_2]=0.$
\item[$L^3_a$] $[x_3,x_1]=x_2,~ [x_3,x_2]=ax_1+x_2,~[x_1,x_2]=0.$
\item[$L^4_a$] $[x_3,x_1]=x_2,~ [x_3,x_2]=ax_1,~[x_1,x_2]=0.$
\end{itemize}
Here the only isomorphisms are $L_a^4\cong L_b^4$ if and only if 
$a=\alpha^2 b$. I have established the non-isomorphism of $L^2$ with
$L_a^3$, $L_a^4$ by Gr\"{o}bner basis calculations. \par
We count the number of non-isomorphic solvable Lie algebras over the
finite field with $q$ elements. The classes $L^1$, $L^2$, $L^3$
always give $q+2$ Lie algebras. If the characteristic of the ground
field is not $2$, then $L_a^4$ gives $3$ more Lie algebras. In that case the
total number is $q+5$. If the characteristic is $2$, then all elements 
of $\mathbb{F}_q$ are squares, meaning that we get two isomorphism classes
of Lie algebras (i.e., $L^4_0$ and $L^4_1$). In that case the total number
is $q+4$.

\begin{rem}
Our classification is the same as the one obtained in \cite{patzas2}.
More precisely, we have $L_{3,1}\cong L^1$, $L_{3,2}\cong 
L_0^4$, $L_{3,3}\cong L_0^3$, $L_{3,4}\cong L^2$, $L_{3,5}\cong
L^4_{\alpha}$ (where $\alpha$ is as in \cite{patzas2}), $L_{3,6}\cong
L^3_{-\alpha}$, and $L_{3,7}\cong L^3_{-1/4}$ (if the characteristic is
not $2$) and $L_{3,7}\cong L^4_1$ (if the characteristic is $2$).\par
So we have the same classification, but with a shorter description.
\end{rem} 

\begin{rem}
From the method used we get an easy algorithm for recognising a given
$3$-dimensional Lie algebra $K$ as one of the $L^i$. First we find a
$2$-dimensional Abelian ideal. Let $x$ span a complement to this ideal. 
Then we find the rational 
canonical form of the adjoint action of $x$ on the ideal. From this we
immediately see to which algebra $K$ is isomorphic.
\end{rem}

\section{The $4$-dimensional case}\label{sec4}

Here we have to find derivation algebras, automorphism groups of
$3$-dimensional Lie algebras $K$. The algebras that will appear in the
final classification will be denoted $M^i$. \par
First we consider $K=L^1$. This Lie algebra is Abelian, so the orbits
of the derivations under the action of $\Aut(K)$ are given by the rational 
canonical form of matrices. If this form consists of three $1\times 1$
blocks, then because of the divisibility condition in Theorem 
\ref{canformthm}, they have to be
the same. After division we get two algebras: the $4$-dimensional
commutative algebra (denoted by $M^1$), and
$$M^2:~~~~~~~~~ 
[x_4,x_1]=x_1, ~ [x_4,x_2]=x_2, ~ [x_4,x_3]=x_3.$$
If there is a $1\times 1$-block and a $2\times 2$-block, then again
because of divisibility we have 
$$D = \begin{pmatrix} s & 0 & 0 \\
                      0 & 0 & -st \\
                      0 & 1 & s+t\\
\end{pmatrix}.$$
Denote the corresponding Lie algebra by $K_{s,t}$. After multiplying
$x_4, x_3$ by $\alpha$ (and $x_1,x_2$ by $1$) we see that this Lie algebra
is isomorphic to $K_{\alpha s,\alpha t}$, where $\alpha\neq 0$. We consider
a few cases. \par
First $s\neq 0$. Then we can take $\alpha= s^{-1}$, and 
we get the Lie algebras 
$$M^3_a:~~~~~~~~~
[x_4,x_1]=x_1, ~ [x_4,x_2]=x_3, ~ [x_4,x_3]=-ax_2+(a+1)x_3.$$
Gr\"{o}bner basis computations reveal that $M^3_a\cong M^3_b$ if and only if
$a=b$. \par
If $s=0$, then if $t\neq 0$, we can take $\alpha = t^{-1}$, and get 
$$M^4:~~~~~~~~~
[x_4,x_2]=x_3, ~ [x_4,x_3]= x_3.$$
Finally, if $s=t=0$ we get
$$M^5:~~~~~~~~~
[x_4,x_2]=x_3.$$
If there is a $3\times 3$-block in the rational normal form, then we get the 
Lie algebras $K_{s,t,u}$:
$$[x_4,x_1] = x_2, ~ [x_4,x_2]=x_3, ~ [x_4,x_3] = sx_1+tx_2+ux_3.$$
Multiplying $x_2,x_3,x_4$ by $\alpha$, $\alpha^2$, $\alpha$ respectively,
we see that $K_{s,t,u}\cong K_{\alpha^3s, \alpha^2t, \alpha u}$.
Hence, if $u\neq 0$, then we can take $\alpha = u^{-1}$, and we get 
the Lie algebras 
$$M^6_{a,b}:~~~~~~~~
[x_4,x_1] = x_2, ~ [x_4,x_2]=x_3, ~ [x_4,x_3] = ax_1+bx_2+x_3.$$
A Gr\"{o}bner basis computation shows that $M^6_{a,b}\cong M^6_{c,d}$
if and only if $a=c$ and $b=d$.\par
If $u=0$ then we get the Lie algebras
$$M^7_{a,b}:~~~~~~~~~
[x_4,x_1] = x_2, ~ [x_4,x_2]=x_3, ~ [x_4,x_3] = ax_1+bx_2.$$
From the above discussion we see that $M^7_{a,b}\cong M^7_{c,d}$ if
$a=\alpha^3 c$ and $b=\alpha^2 d$ (for some $\alpha\neq 0$). From a 
Gr\"{o}bner basis computation we have that this is also a necessary
condition. So, if both parameters are nonzero, than by a suitable choice
for $\alpha$ we can make them equal. Hence this class splits into three
subclasses: $M^7_{a,a}$, $M^7_{a,0}$, $M^7_{0,b}$. Among the first class 
there are no isomorphisms. \par
Now we consider $K=L^2$. The coset representatives of the outer derivations
of $K$ (modulo inner derivations) are
$$D= \begin{pmatrix} s & t & 0 \\ 
                     u & 0 & 0 \\
                     0 & 0 & 0 \\
     \end{pmatrix}.$$
The automorphism group of $K$ consists of the elements 
$$\phi = \begin{pmatrix} \alpha & \beta & \epsilon_1 \\  
                         \gamma & \delta & \epsilon_2 \\
                          0 & 0 & 1 \\
         \end{pmatrix},$$
where $\alpha\delta - \beta\gamma \neq 0$. Modulo inner derivations we
have
$$(\det \phi) \phi D\phi^{-1} = \begin{pmatrix}
  s\alpha\delta - t\alpha\gamma + u\beta\delta  & 
  -s\alpha\beta + t\alpha^2 - u\beta^2 & 0 \\
 s\gamma\delta - t\gamma^2 + u\delta^2  &
 -s\beta\gamma + t\alpha\gamma - u\beta\delta & 0 \\
0&0&0\\
\end{pmatrix}.$$       
First we suppose that $s\neq 0$. Then we choose $\beta=\gamma=0$ and
$\alpha =1/s$, $\delta=1$. Then $D$ is conjugate to 
$$\begin{pmatrix} 1 & w & 0\\
                  v & 0 & 0\\
                  0 & 0 & 0 \\
\end{pmatrix}$$
for some $v,w$. This leads to the Lie algebras $K_{v,w}$ with basis
$x_1,x_2,x_3,x_4$ and nonzero commutators
$$[x_4,x_1] = x_1+vx_2, ~ [x_4,x_2]=wx_1,~ [x_3,x_1]=x_1,~ [x_3,x_2]=x_2.$$
Suppose that $w\neq 0$, then by setting $y_1= wx_1$, $y_i = x_i$ for
$i=2,3,4$, we see that $K_{v,w}\cong K_{v',1}$. Denote this Lie algebra
simply by $K_v$. By some calculations
it is seen that the centralizer $C(\ad K_v)$ in the full (associative)
matrix algebra is spanned by the identity and 
$$\begin{pmatrix} 1 & 1 & 0 & 0 \\
                  v & 0 & 0 & 0 \\
                  0 & 0 & 0 & v \\
                  0 & 0 & 1 & 1 \\
 \end{pmatrix}.$$
The minimal polynomial of this last matrix is $T^2-T-v$. Suppose
that the characteristic of $F$ is not $2$. Then, if 
$v=-\frac{1}{4}$, this algebra has a nonzero radical. We get the Lie
algebra
$$N:~~~~~~~~~
[x_4,x_1] = x_1-\frac{1}{4}x_2, ~ [x_4,x_2]=x_1,~ [x_3,x_1]=x_1,~ 
[x_3,x_2]=x_2.$$
(We denote this algebra by $N$ and not by $M^8$, because it is isomorphic
to a Lie algebra that we define later).
On the other hand, if $v\neq -\frac{1}{4}$ then $C(\ad K_v)$ is semisimple. 
Also, if $T^2-T-v$ has a root in the base field, then it splits. This implies
that $K_v$ is isomorphic to the direct sum of two $2$-dimensional Lie
algebras (namely the non-commutative ones) (cf. \cite{gra6}, \S 1.15, 
\cite{rwz}). We get the Lie algebra
$$M^8:~~~~~~~~~
[x_1,x_2]=x_2,~ [x_3,x_4]=x_4.$$
Now suppose that $T^2-T-v$ does not have a root in $F$. Then
$K_v$ is indecomposable. Suppose that $K_v\cong K_w$, where also 
$T^2-T-w$ has no root in $F$. Then from the
Gr\"{o}bner basis it follows that $v+\frac{1}{4}= \alpha^2(w+
\frac{1}{4})$ for some nonzero $\alpha\in F$.
(There is also another argument to prove this: 
as seen above $K_v$ splits over $F(\sqrt{1+4v})$ so also $K_w$ splits 
over this field. Hence $\sqrt{1+4w}\in F(\sqrt{1+4v})$. This implies the
claim.) Conversely, suppose that $v+\frac{1}{4}= \alpha^2(w+\frac{1}{4})$ 
for some nonzero $\alpha\in F$. Let $\phi: K_v\to K_w$ be the linear map 
given by 
$\phi(x_1) =\alpha y_1 +\frac{1}{2}(1-\alpha)y_2$, $\phi(y_2) = y_2$, 
$\phi(y_3) = y_3$, $\phi(x_4) = \frac{1}{2}(1-\alpha)y_3+\alpha y_4$. Then 
$\phi$ is an isomorphism. We conclude that $K_v\cong K_w$ if and only
if $v+\frac{1}{4}= \alpha^2(w+\frac{1}{4})$ for some nonzero 
$\alpha\in F$. \par
Now we deal with the case where the characteristic of $F$ is $2$. 
Just as above, if $T^2+T+a$ factors over $F$, then $K_v$ is isomorphic
to a direct sum. If the polynomial does not factor, then $K_v$ is
indecomposable. From the Gr\"{o}bner basis computation it follows that
$K_v\cong K_w$ implies that $X^2+X+v+w$ has roots in $F$. Conversely,
suppose that this equation has a root $\alpha\in F$. Then there is an
isomorphism $\phi : K_v\to K_w$ given by $\phi(x_1)=y_1+\alpha y_2$,
$\phi(x_2)=y_2$, $\phi(x_3)=y_3$, $\phi(x_4) = \alpha y_3+y_4$. So
$K_v\cong K_w$ if and only if $X^2+X+v+w$ has roots in $F$. The 
conclusion is that we get the Lie algebras
$$M^9_a:~~~~~~~~~
[x_4,x_1] = x_1+ax_2, ~ [x_4,x_2]=x_1,~ [x_3,x_1]=x_1,~ [x_3,x_2]=x_2,$$
where $a\in F$ is such that $T^2-T-a$ has no root in the base field.\par
Now suppose that $w=0$. Then $K_{v,0}$ is the direct sum of ideals with
bases $x_1+vx_2$, $x_4$, and $x_2,x_3-x_4$. So $K_{v,0}\cong M^8$. This 
finishes the case where $s\neq 0$.\par
Now suppose that $s=0$.
Then $D$ is equal to 
$$\begin{pmatrix} 0 & t & 0\\
                  u & 0 & 0\\
                  0 & 0 & 0 \\
\end{pmatrix}$$
If $u\neq 0$ then we can divide by it and obtain the derivation
$$\begin{pmatrix} 0 & a & 0\\
                  1 & 0 & 0\\
                  0 & 0 & 0 \\
\end{pmatrix}$$
This leads to the Lie algebras 
$$M_a^{10}:~~~~~~~~~
[x_4,x_1] = x_2,~ [x_4,x_2]=ax_1,~ [x_3,x_1]=x_1,~ [x_3,x_2]=x_2.$$
If the characteristic of $F$ is not $2$, then 
$M_a^{10}\cong M^9_{a-\frac{1}{4}}$. The isomorphism is given by 
$\phi(x_1)=2y_2$, $\phi(x_2)=2y_1-y_2$, $\phi(x_3)=y_3$, 
$\phi(x_4)=-\frac{1}{2}y_3+y_4$ (where the $x_i$ are the basis
elements of $M_a^{10}$). Note that, if $a=0$, this gives an isomorphism
with $N$.\par
If the characteristic is $2$, then from a Gr\"{o}bner basis computation
it follows that $M_a^{10}$ is not isomorphic to $M^9_b$. So in this case we
have a new series of Lie algebras.
From a Gr\"{o}bner basis computation we get that $M^{10}_a\cong 
M^{10}_b$ implies that $Y^2+X^2b+a=0$ is solvable in $F$, with
$X\neq 0$. On the other hand, if $\beta$ and $\alpha\neq 0$ are such
that $\beta^2+\alpha^2b+a=0$, then $\phi(x_1) = y_1$, $\phi(x_2) = 
\beta y_1+\alpha y_2$, $\phi(x_3) = y_3$, $\phi(x_4) =\beta y_3+\alpha y_4$
is an isomorphism. So $M^{10}_a\cong M^{10}_b$ if and only if
$Y^2+X^2b+a=0$ has a solution in $F$, with $X\neq 0$. In particular, if 
the field is perfect (i.e., $F^2=F$) then $M_a^{10}\cong M_0^{10}$. \par
If $u=0$, and $t\neq 0$ then we divide by $t$. The corresponding Lie algebra
has multiplication table
$$[x_4,x_2]=x_1,~ [x_3,x_1]=x_1,~ [x_3,x_2]=x_2.$$
If the characteristic is not $2$, then it is isomorphic to $N$, the
isomorphism being $\phi(x_1)=2y_1-y_2$, $\phi(x_2)=y_2$, $\phi(x_3)=y_3$,
$\phi(x_4)=-y_3+2y_4$. If the characteristic is $2$, then this algebra is
isomorphic to $M^{10}_0$, with in this case, $\phi(x_1)=y_2$,
$\phi(x_2)=y_1$, $\phi(x_3)=y_3$, $\phi(x_4)=y_4$.\par
If $t=0$ as well, then the derivation is inner, and we obtain nothing 
new.\par
Now we consider the Lie algebra $K=L^3_a$. Its derivations consist of
$$\begin{pmatrix} u & av & s\\
                  v & u+v & t\\
                  0 & 0 & 0 \\
\end{pmatrix}.$$
This means that, if $a\neq 0$, then modulo scalar factors there is only
one outer derivation, namely
$$D=\begin{pmatrix} 1 & 0 & 0\\
                  0 & 1 & 0\\
                  0 & 0 & 0 \\
\end{pmatrix},$$
leading to the Lie algebras
$$[x_4,x_1] = x_1,~ [x_4,x_2]=x_2,~ [x_3,x_1]=x_2,~ [x_3,x_2]=ax_1+x_2.$$
However, by interchanging $x_3,x_4$ and $x_1,x_2$ we get the Lie algebra
$K_a$ considered before (leading to the algebras $N$, $M^8$, $M^9_a$).
\par
If $a=0$, then apart from the derivation above, we get two more:
$$D_1 = \begin{pmatrix} 1 & 0 & v\\
                  0 & 1 & 0\\
                  0 & 0 & 0 \\
\end{pmatrix}, \text{~~ and ~~}
D_2 = \begin{pmatrix} 0 & 0 & 1\\
                  0 & 0 & 0\\
                  0 & 0 & 0 \\
\end{pmatrix}.$$
Let $\sigma$ be the automorphism with matrix
$$\begin{pmatrix} 1 & 0 & v\\
                  0 & 1 & 0\\
                  0 & 0 & 1 \\
\end{pmatrix}.$$
Then $\sigma D_1\sigma^{-1}$ is equal to $D$, so we get nothing new from
$D_1$. However, we can not get rid of $D_2$ in this way. It leads to the
Lie algebra
$$[x_4,x_3] = x_1,~ [x_3,x_1]=x_2,~ [x_3,x_2]=x_2.$$
But this is isomorphic to $M^6_{0,0}$, by $\phi(x_1)=y_2$, $\phi(x_2)=y_3$,
$\phi(x_3)=y_4$, $\phi(x_4)=-y_1$.\par
Now we deal with the Lie algebra $K=L_a^4$. If $a\neq 0$ and 
the characteristic of $F$ is not $2$, then the derivations of 
$L_a^4$ are given by
$$\begin{pmatrix} u & av & s\\
                  v & u  & t\\
                  0 & 0 & 0 \\
\end{pmatrix}.$$
Here modulo inner derivations, and scalar factors, there remains only
one derivation 
$$\begin{pmatrix} 1 & 0 & 0\\
                  0 & 1 & 0\\
                  0 & 0 & 0 \\
\end{pmatrix},$$
leading to the Lie algebras
$$[x_4,x_1] = x_1,~ [x_4,x_2]=x_2,~ [x_3,x_1]=x_2,~ [x_3,x_2]=ax_1.$$
This Lie algebra is isomorphic to $M^9_{a-\frac{1}{4}}$. The isomorphism
is given by $\phi(x_1)=2y_2$, $\phi(x_2)=2y_1-y_2$, $\phi(x_3)=
-\frac{1}{2}y_3+y_4$, $\phi(x_4)=y_3$.\par
If $a\neq 0$ and the characteristic of $F$ is $2$, then the derivations
are given by 
$$\begin{pmatrix} u & av & s\\
                  v & u+w & t\\
                  0 & 0 & w \\
\end{pmatrix}.$$
So modulo inner derivations we get 
$$\begin{pmatrix} u & 0 & 0\\
                  0 & u+w  & 0\\
                  0 & 0 & w \\
\end{pmatrix}.$$
If $w=0$ then this leads to the algebra that we have seen in the 
characteristic not $2$ case. In this case it is isomorphic to 
$M^{10}_a$, by $\phi(x_1)=y_1$, $\phi(x_2)=y_2$, $\phi(x_3)=y_4$,
$\phi(x_4)=y_3$.\par
If $u$ and $w$ are both nonzero, then after dividing we get
the derivations
$$\begin{pmatrix} 1 & 0 & 0\\
                  0 & b & 0\\
                  0 & 0 & 1+b \\
\end{pmatrix}.$$
Here we assume that $b\neq 1$ as we have already listed the corresponding
algebra (it is isomorphic to $M^{10}_a$).
The Lie algebras we now get are:
$$M^{11}_{a,b}:~~~~~~~~~
[x_4,x_1] = x_1, [x_4,x_2] = bx_2,~ [x_4,x_3]=(1+b)x_3,~ [x_3,x_1]=x_2,~ 
[x_3,x_2]=ax_1.$$
(Here $a\neq 0$, $b\neq 1$.) Let $c\neq 0$ and $d\neq 1$.
Set $\delta = (b+1)/(d+1)$. We claim that 
$M^{11}_{a,b}\cong M^{11}_{c,d}$ if and only if $(\delta^2+(b+1)\delta+b)/c$ 
and
$a/c$ are squares in $F$. The only if part follows from inspection of the 
Gr\"{o}bner basis. The if part from explicit construction of an isomorphism. 
Let $\gamma,\epsilon\in F$ be such that
$$\gamma^2 = \frac{1}{c}(\delta^2+(b+1)\delta+b), \text{ and } 
\epsilon^2 = \frac{a}{c}.$$
If $\delta =1$ then $b=d$ and isomorphism follows already from the 
dimension $3$ isomorphism. So we suppose that $\delta\neq 1$, and we set
$\beta=\delta+1$, $\alpha = c\gamma$. Then $\phi : K_{a,b} \to K_{c,d}$
given by $\phi(x_1) =\alpha y_1+\beta y_2$, $\phi(x_2) = c\epsilon\beta y_1
+\alpha\epsilon y_2$, $\phi(x_3)=\epsilon y_3$, $\phi(x_4) = 
\gamma y_3+\delta y_4$, is an isomorphism. In particular, if $F$ is perfect,
then $M^{11}_{a,b}\cong M^{11}_{1,0}$. \par
If $w\neq 0$, but $u=0$, then we get the algebra
$$[x_4,x_2] = x_2,~ [x_4,x_3]=x_3,~ [x_3,x_1]=x_2,~ [x_3,x_2]=ax_1.$$
If $a\neq 0$ and $a\neq 1$ then this is isomorphic to $M^{11}_{a,a}$,
by $\phi(x_1)=y_1+a^{-1}y_2$, $\phi(x_2)=y_1+y_2$, $\phi(x_3) = 
y_3$, $\phi(x_4)=\frac{1}{a+1}(y_3+y_4)$. If $a=1$, then it is isomorphic
to $M^{11}_{1,0}$, by $\phi(x_1)=y_2$, $\phi(x_2)=y_1$, $\phi(x_3) = 
y_3$, $\phi(x_4)=y_4$. If $a=0$, then it is isomorphic to $M^{10}_0$,
by $\phi(x_1)=y_4$, $\phi(x_2)=y_2$, $\phi(x_3)=y_1$, $\phi(x_4)=y_3$.\par
Now suppose that $a=0$. Then the derivations (modulo inner derivations) are 
$$D= \begin{pmatrix} u_1 & 0 & v_1\\
                  0 & u_1+v_3 & 0\\
                  u_3 & 0 & v_3 \\
\end{pmatrix}.$$
A general automorphism of $L_0^4$ is given by 
$$\phi = \begin{pmatrix}
\alpha_1 & 0 & \gamma_1 \\
\alpha_2 & \alpha_1\gamma_3-\alpha_3\gamma_1 & \gamma_2 \\
\alpha_3 & 0 & \gamma_3 \\
\end{pmatrix},$$
where $\alpha_1\gamma_3-\alpha_3\gamma_1\neq 0$.
Now the entry at position $(3,3)$ of $(\alpha_1\gamma_3-\alpha_3\gamma_1)
\phi D \phi^{-1}$ is $-u_1\alpha_3\gamma_1 -u_3\gamma_1\gamma_3 +
v_1\alpha_1\alpha_3 + v_3\alpha_1\gamma_3$. It is straightforward to see 
that, except in the case where $u_1=v_3=s\neq 0$ and $u_3=v_1=0$,
we can choose the $\alpha_i,\gamma_i$ such that this becomes zero.
In the former case we divide by $s$ and get the Lie algebra
$$M^{12}:~~~~~~~~~
[x_4,x_1] = x_1,~ [x_4,x_2]=2x_2, ~ [x_4,x_3] = x_3,~ [x_3,x_1]=x_2.$$
Otherwise $D$ is conjugate to 
$$D' = \begin{pmatrix} u_1 & 0 & v_1\\
                  0 & u_1 & 0\\
                  u_3 & 0 & 0 \\
\end{pmatrix}.$$ 
Now let $\phi$ be the automorphism given by the matrix
$$\begin{pmatrix} \alpha & 0 & 0\\
                  0 & \alpha\beta & 0\\
                  0 & 0 & \beta \\
\end{pmatrix},$$
where both $\alpha,\beta$ are nonzero. Then
$$\phi D'\phi^{-1} = \begin{pmatrix} u_1 & 0 & \frac{\alpha}{\beta}v_1\\
                  0 & u_1 & 0\\
                  \frac{\beta}{\alpha} u_3 & 0 & 0 \\
\end{pmatrix}.$$
First suppose that $u_1\neq 0$. Then we divide by it. If $v_1$ is also
nonzero, then we choose $\alpha = u_1$ and $\beta=v_1$. This leads to the
Lie algebras
$$M^{13}_a:~~~~~~~~~
[x_4,x_1] = x_1+ax_3,~ [x_4,x_2]=x_2, ~ [x_4,x_3] = x_1,~ [x_3,x_1]=x_2.$$
From the Gr\"{o}bner basis it follows that two of those algebras, 
with parameters $a$ and $b$, are isomorphic if and only if $a=b$.\par
Now suppose that $v_1=0$. If $u_3\neq 0$, then set $\alpha = u_3$, 
$\beta=u_1$. We get the Lie algebra
$$[x_4,x_1] = x_1+x_3,~ [x_4,x_2]=x_2, ~ [x_3,x_1]=x_2.$$
If we set $\tilde{x}_1 = x_1+x_3$, $\tilde{x}_2=-x_2$, $\tilde{x}_3 = 
x_1$, $\tilde{x}_4 = x_4$, then we see that with respect to this new basis
the Lie algebra has the same multiplication table as $M^{13}_0$.\par
Suppose that $u_3$ is zero as well. Then we get the Lie algebra
$$[x_4,x_1] = x_1,~ [x_4,x_2]=x_2, ~ [x_3,x_1]=x_2.$$
In this case we set $\tilde{x}_1 = x_1$, $\tilde{x}_2=x_2$, $\tilde{x}_3 = 
x_1+x_3$, $\tilde{x}_4 = x_4$. Again we get the multiplication table of
$M^{13}_0$.\par
Secondly, suppose that $u_1=0$. If $v_1\neq 0$, then we divide by it,
and get the derivations
$$\begin{pmatrix} 0 & 0 & 1\\
                  0 & 0 & 0\\
                  a & 0 & 0 \\
\end{pmatrix},$$
leading to the Lie algebras
$$M^{14}_a:~~~~~~~~~
[x_4,x_1] = ax_3,~ [x_4,x_3]=x_1, ~ [x_3,x_1]=x_2.$$
By setting
$y_i= \alpha x_i$ for $i=1,2,4$, $y_3=x_3$, we see that this Lie algebra
is isomorphic to the same one with parameter $\alpha^2a$. On the other hand,
from a Gr\"{o}bner basis computation we get that $M^{14}_a\cong M^{14}_b$
implies $a= \alpha^2b$ for some $\alpha$.\par
If $v_1=0$, then we get two more algebras. The first is 
a direct sum isomorphic to $M^5$. The other is 
$$[x_4,x_1] = x_3, ~ [x_3,x_1]=x_2.$$
Here we set $\tilde{x}_1 = x_3$, $\tilde{x}_2=-x_2$, $\tilde{x}_3 = 
x_1$, $\tilde{x}_4 = x_4$. This gives us the multiplication table of
$M^{14}_0$.\par
There are some additional isomorphisms. If the characteristic of $F$ 
is not $2$, then $M^{13}_0\cong N$, by
$\phi(x_1)=y_1+y_2$, $\phi(x_2)=3y_1-\frac{3}{2}y_2$, $\phi(x_3) =
y_1+y_2-y_3+2y_4$, $\phi(x_4)=y_3$. For this reason we don't list
$N$ separately (if the characteristic is not $2$, then it doesn't exist).\par
Second, $M^7_{0,0}\cong M^{14}_0$, by $\phi(x_1)=-y_4$, $\phi(x_2)=y_1$, 
$\phi(x_3)=y_2$, $\phi(x_4)=y_3$.\par
If the characteristic is $2$, and $a=\alpha^2\neq 0$, then 
$M_a^{10}\cong M^{13}_0$, by $\phi(x_1)=y_1$, $\phi(x_2)=\alpha y_1+
\alpha y_2$, $\phi(x_3)=y_1+y_2+y_4$, $\phi(x_4)=\alpha y_3+
\alpha y_4$. Note that if $a=0$ we also have isomorphism: 
$M_0^{10}\cong M_1^{10}\cong M_0^{13}$.\par
I have established the non-isomorphism of the remaining Lie algebras
$M^i$ by Gr\"{o}bner basis computations.

\begin{rem}
The Lie algebras $M^{11}_{0,b}$ do exist, and one may wonder where
they occur in the list. We have $M^{11}_{0,0}\cong M^{12}$, and
$M^{11}_{0,b}\cong M^{13}_{(b+1)/b^2}$ if $b\neq 0, b\neq 1$. 
\end{rem}

We count the number of solvable Lie algebras of dimension $4$ over the finite
field $\F_q$, where $q=p^m$. For that we start with a well-known lemma 
(see \cite{ber}, Theorems 6.69, 6.695).

\begin{lem}\label{lem3}
Let $u\in \F_q$, where $q=2^m$. Then the equation $X^2+X+u$ has a solution 
in $\F_q$ if and only if 
$$\Tr_2(u) = \sum_{i=0}^{m-1} u^{2^i}=0.$$
\end{lem}

The classes $M^1$, $M^2$, $M^3_a$, $M^4$, 
$M^5$, $M^6_{a,b}$ contain 1, 1, $q$, 1, 1, $q^2$ algebras respectively.
As noted before, the class $M^7_{a,b}$ splits in three subclasses:
$M^7_{a,a}$, $M^7_{a,0}$ ($a\neq 0$), and $M^7_{0,b}$ ($b\neq 0$). The first 
of these contains $q$ elements. We have $M^7_{a,0}\cong M^7_{a',0}$ if and
only if $a=\alpha^3 a'$ for some $\alpha\in \F_q$. First suppose that $q$ is odd.
If $q\equiv 1\bmod 6$, then $X^3=1$ has $3$ solutions in $\F_q$. In that case 
$\F_q$ contains $(q-1)/3$ cubes and hence we get $3$ algebras. If $q \not\equiv 
1\bmod 6$ then $X^3=1$ has $1$ solution in $\F_q$ and hence $\F_q^3=\F_q$ 
and we get only $1$ algebra. Now suppose that $p=2$, then 
$\Tr_2(1) = m$. So by Lemma \ref{lem3}, $X^2+X+1$ has solutions in $\F_q$
if and only if $m$ is even. This is the same as saying that $p^m \equiv 4 \bmod 6$.
In the same way as above we conclude that in this case $M^7_{a,0}$ has $3$
algebras. In the other case, $q \equiv 2\bmod 6$, we get $1$ algebra. 
The class $M^7_{0,b}$ contains $1$ algebra if $p=2$, and $2$ algebras
if $p>2$. Summarising
$$|M^7_{a,b}| = \begin{cases}
                 q+5 & q\equiv 1 \bmod 6\\
                 q+2 & q\equiv 2 \bmod 6\\
                 q+3 & q\equiv 3 \bmod 6\\
                 q+4 & q\equiv 4 \bmod 6\\
                 q+3 & q\equiv 5 \bmod 6
                \end{cases}. $$
From $M^8$ we get $1$ algebra. Now we consider $M^9_a$. First suppose that 
$q$ is odd. We have to find the set of $a\in \F_q$ such that $T^2-T-a$ has no 
root in $\F_q$. Suppose that this equation has a root $\alpha$. Then the other 
root is $1-\alpha$ and $a=\alpha^2-\alpha$. Let $B$ be the set of all 
$\alpha^2-\alpha$ for $\alpha \in \F_q$. If the equation $X^2-X=c$ has one 
solution in $\F_q$, then it has two solutions, unless $c=-\frac{1}{4}$.
This implies that $|B|=(q+1)/2$. Let $A$ be the set of all $a\in \F_q$
such that $T^2-T-a$ has no root in $\F_q$. Then $|A|=(q-1)/2$. Also, for
$0\neq \beta\in\F_q$ we define $h_{\beta} : \F_q\to \F_q$ by $h_{\beta}(x) = 
\beta^2(x+\frac{1}{4})-\frac{1}{4}$. Then $h_{\beta}$ is a bijection. It
stabilizes $B$ and hence $A$. Now $M^9_a\cong M^9_b$ precisely if $h_{\beta}(a)
=b$ for some $0\neq \beta\in\F_q$. There are exactly $(q-1)/2$ different
$h_{\beta}$'s. So all $M^9_a$ for $a\in A$ are isomorphic. Hence we get $1$ 
algebra. Secondly, suppose that $q$ is even. Choose $a,b$ such that $T^2+T+a$
and $T^2+T+b$ have no roots in $\F_q$. Then by Lemma \ref{lem3},
$\Tr_2(a)=\Tr_2(b)=1$. But then $\Tr_2(a+b)=0$ and $T^2+T+a+b$ has roots
in $\F_q$. Hence $M^9_a\cong M^9_b$. So we get one algebra in this case as well.\par
The classes $M^{10}_a$ and $M^{11}_{a,b}$ are only defined for characteristic
$2$. For perfect fields they both have one algebra. However, the algebra 
in $M^{10}_a$ disappears due to the isomorphism with $M^{13}_0$. \par
The classes $M^{12}$ and $M^{13}_a$ have $1$ and $q$ algebras respectively. 
For $M_a^{14}$ we exclude $a=0$, as that algebra is isomorphic to $M^7_{0,0}$.
Therefore this class contains $1$ algebra if $p=2$ and $2$ algebras if 
$p>2$. \par
Now we add these numbers, and find that the total number of solvable Lie algebras 
over $\F_q$ is
$$q^2+3q+9+    \begin{cases}
                 5 & q\equiv 1 \bmod 6\\
                 2 & q\equiv 2 \bmod 6\\
                 3 & q\equiv 3 \bmod 6\\
                 4 & q\equiv 4 \bmod 6\\
                 3 & q\equiv 5 \bmod 6
                \end{cases}, $$
which is slightly more than the number found in \cite{patzas2}.

\begin{rem}
With
$L_{4,i}$ as in \cite{patzas2} we have $L_{4,1}\cong M^1$, $L_{4,2}\cong
M_5$ (for this one has to correct the table given in \cite{patzas2}; with 
the table as given in \cite{patzas2}, we have $L_{4,2}\cong L_{4,3}$), 
$L_{4,3}\cong M^7_{0,0}$, $L_{4,4}\cong M^4$, $L_{4,5}\cong
M_0^3$, $L_{4,6}\cong M^7_{0,\alpha}$, $L_{4,7}\cong M^6_{0,-\alpha}$, 
$L_{4,8}\cong M^8$, $L_{4,9}\cong M^6_{0,-\frac{1}{4}}$ (if the
characteristic is not $2$), $L_{4,9}\cong M^7_{0,1}$ (if the characteristic
is $2$), $L_{4,10}\cong M^2$, $L_{4,11}\cong M_{\alpha}^3$, 
$L_{4,12}\cong M^6_{\alpha_3,-\alpha_2}$, $L_{4,13}\cong M^7_{\alpha,\alpha}$,
$L_{4,14}\cong M^7_{\alpha,0}$, $L'_{4,8}\cong M^9_{-\alpha}$,
$L_{4,15}\cong M^7_{-2,3}$ (characteristic not $3$) $L_{4,15}\cong
M_1^3$ (characteristic $3$), $L_{4,16}\cong M^6_{-2\alpha^3+\alpha^2,
3\alpha^2-2\alpha}$ ($\alpha\neq \frac{1}{3}$), $L_{4,16}\cong M^3_1$
($\alpha = \frac{1}{3}$), $L_{4,17}\cong M^6_{\frac{1}{27},-\frac{1}{3}}$
(characteristic not $3$ and $\alpha\neq 0$), $L_{4,17}\cong M_{1,0}^7$ 
(characteristic $3$, and $\alpha\neq 0$), $L_{4,17}\cong M^7_{0,0}$
($\alpha=0$), $L_{4,18}\cong M^{12}$, $L_{4,19}\cong M_{\alpha}^{14}$,
$L_{4,20}\cong M^{13}_{-\alpha}$. \par
In \cite{patzas2} the algebra $M^{11}_{a,b}$ is missing.
This can be explained by the circumstance that the
method used in \cite{patzas2} relies on the derived algebra being
nilpotent. Now, if $a\neq 0$ and $b\neq 1$ then the derived algebra
of $M^{11}_{a,b}$ is not nilpotent.
\end{rem}

\begin{rem}
Of course the next step will be to describe the classification for
dimension $5$. This will be the subject of a next paper. A problem
that may occur is that the computation of the coordinates of the 
elements of a Gr\"{o}bner basis (with respect to the input basis)
makes the algorithm rather slow. So there may occur cases where it is
impossible to compute a Gr\"{o}bner basis (with coordinates).
\end{rem}

\begin{rem}
As in the dimension $3$ case, it is possible to formulate an algorithm
that for a given solvable Lie algebra $K$ of dimension $4$ finds the $M^i$ to which
it is isomorphic. First we find a $3$-dimensional ideal, and establish
to which of the $L^i$ it is isomorphic. Then for each of the four
possibilities we basically follow the classification procedure. 
\end{rem}

\bibliographystyle{plain}
\def\cprime{$'$} \def\cprime{$'$} \def\cprime{$'$}

\end{document}